\pdfoutput=1
\RequirePackage{ifpdf}
\ifpdf 
\documentclass[pdftex]{sigma}
\else
\documentclass{sigma}
\fi

\newcommand{\HVP}[1]{\widetilde{X}_{#1}}
\newcommand{\CHVP}[1]{\widetilde{\alpha}_{#1}}
\newtheorem{Th}{Theorem}

\newtheorem{Proposition}[Th]{Proposition}
\newtheorem{Lemma}[Th]{Lemma}
{\theoremstyle{definition}
\newtheorem{Rem}[Th]{Remark}}
\DeclareMathOperator{\Lie}{Lie}

\begin{document}

\allowdisplaybreaks

\renewcommand{\PaperNumber}{048}

\FirstPageHeading

\ShortArticleName{Multi-Hamiltonian Structures on Spaces of Closed Equicentroaf\/f\/ine Plane Curves}

\ArticleName{Multi-Hamiltonian Structures on Spaces\\
of Closed Equicentroaf\/f\/ine Plane Curves\\
Associated to Higher KdV Flows}

\Author{Atsushi FUJIOKA~$^\dag$ and Takashi KUROSE~$^\ddag$}

\AuthorNameForHeading{A.~Fujioka and T.~Kurose}

\Address{$^\dag$~Department of Mathematics, Kansai University, Suita, 564-8680, Japan}
\EmailD{\href{mailto:afujioka@kansai-u.ac.jp}{afujioka@kansai-u.ac.jp}}

\Address{$^\ddag$~Department of Mathematical Sciences, Kwansei Gakuin University, Sanda, 669-1337, Japan}
\EmailD{\href{mailto:crg31562@kwansei.ac.jp}{crg31562@kwansei.ac.jp}}

\ArticleDates{Received October 11, 2013, in f\/inal form April 16, 2014; Published online April 22, 2014}

\Abstract{Higher KdV f\/lows on spaces of closed equicentroaf\/f\/ine plane curves are studied and it is shown that the f\/lows
are described as certain multi-Hamiltonian systems on the spaces.
Multi-Hamiltonian systems describing higher mKdV f\/lows are also given on spaces of closed Euclidean plane curves via the
geometric Miura transformation.}


\Classification{37K25; 35Q53}

\section{Introduction}

A motion of a~curve is a~smooth one-parameter family of connected curves in a~space.
It is known that many dif\/ferential equations related to integrable systems can be linked with special motions of
curves~\cite{CQ1,CQ2,CQ3,RS}.
For example, for a~special motion of an inextensible curve in the Euclidean plane, the curvature evolves according to
the modif\/ied Korteweg--de Vries (mKdV) equation~\cite{Lam} (cf.~Section~\ref{Section4} below).
There are a~lot of preceding studies on motions of curves related to Euclidean geometry and the mKdV equation.
See~\cite{BSW,SW,TT} and references therein.
For special motions of a~space curve, it is also known that the nonlinear Schr\"odinger equation appears~\cite{H}.
In~\cite{FK1,FK2}, the authors studied motions of a~curve in the complex hyperbola under which the curvature evolves
according to the Burgers equation.

In this paper, we shall study motions of an equicentroaf\/f\/ine plane curve.
Under a~special motion of an equicentroaf\/f\/ine plane curve, the equicentroaf\/f\/ine curvature evolves according to the
Korteweg--de Vries (KdV) equation.
In order to explain the above motion geometrically, Pinkall~\cite{P} introduced the natural presymplectic form on the
space of closed equicentroaf\/f\/ine plane curves with f\/ixed enclosing area, and showed that the equicentroaf\/f\/ine curvature
evolves according to the KdV equation when the f\/low is generated by the total equicentroaf\/f\/ine curvature.
Furthermore, the result has been generalized to the case of higher KdV f\/lows (cf.~\cite{CIB,FK3}).

On the other hand, it is known that a~lot of completely integrable systems are described as bi-Hamiltonian systems, from
which the existence of many f\/irst integrals can be deduced (Magri's theorem~\cite{M,O}).
In this context, many of motions of curves as above have been studied from the viewpoint of bi-Hamiltonian systems
recently~\cite{A1,A2,A3,A4,AA1,AA2,AM,AV,LQZ,B,BSW,SB}.
The purpose of this paper is to construct a~multi-Hamiltonian structure associated to the higher KdV f\/lows on each level
set of Hamiltonian functions in a~geometric way (Theorem~\ref{Th:MHSForKdV}).
Moreover, we shall also introduce multi-Hamiltonian structures associated to the higher mKdV f\/lows on the spaces of
closed Euclidean plane curves via the geometric Miura transformation.

\section{A bi-Hamiltonian structure on the space\\of closed equicentroaf\/f\/ine curves}

Throughout this paper all maps are assumed to be smooth.

For a~regular plane curve $\gamma$ whose velocity vector is transversal to the position vector at each point, we can
choose the parameter~$s$ of $\gamma$ as $\det
\begin{pmatrix}
\gamma(s)
\\
\gamma_s(s)
\end{pmatrix}
\equiv 1$ holds.
A plane curve $\gamma$ provided with such a~parameter~$s$ is called an \textit{equicentroaffine} plane curve.
For an equicentroaf\/f\/ine plane curve $\gamma$, we can def\/ine a~function $\kappa$,
called the \textit{equicentroaffine curvature}, by $\gamma_{ss}=-\kappa\gamma$.

We set the space $\mathcal{M}$ of closed equicentroaf\/f\/ine plane curves by
\begin{gather*}
\mathcal{M}=\left\{\gamma:S^1\to\mathbb{R}^2\setminus\{0\}\left\vert\det
\begin{pmatrix}
\gamma
\\
\gamma_s
\end{pmatrix}
=1\right.\right\},
\end{gather*}
where $S^1=\mathbb{R}/2\pi\mathbb{Z}$.
Let $\gamma(\, \cdot\,,t)\in\mathcal{M}$ be a~one-parameter family of closed equicentroaf\/f\/ine plane curves.
As in~\cite{P}, the motion vector f\/ield $\gamma_t$ is represented as
\begin{gather}\label{Eq:2.1}
\gamma_t =-\frac{1}{2}\alpha_s\gamma+\alpha\gamma_s,
\qquad
\alpha:\ S^1\to\mathbb{R},
\end{gather}
and the equicentroaf\/f\/ine curvature $\kappa$ evolves as
\begin{gather}\label{Eq:2.2}
\kappa_t=\Omega\alpha_s=\frac{1}{2}\alpha_{sss}+2\kappa\alpha_s+\kappa_s\alpha,
\end{gather}
where
\begin{gather*}
\Omega=\frac{1}{2}D_s^2+2\kappa+\kappa_sD_s^{-1},
\qquad
D_s = \frac{\partial}{\partial s},
\end{gather*}
is the \textit{recursion operator of the KdV equation}:
\begin{gather*}
\kappa_t=\Omega\kappa_s=\frac{1}{2}\kappa_{sss}+3\kappa\kappa_s.
\end{gather*}
Hence when we choose the one-parameter family $\gamma(\,\cdot, t)$ as $\alpha=D_s^{- 1}\Omega^{n-1}\kappa_s$, we obtain
the~$n$th KdV equation for $\kappa$:
\begin{gather}\label{Eq:2.4}
\kappa_t=\Omega^n\kappa_s.
\end{gather}

The tangent space of $\mathcal{M}$ at $\gamma\in\mathcal{M}$ is described as
\begin{gather*}
T_{\gamma}\mathcal{M}=\left\{\left.{-}\frac{1}{2}\alpha_s\gamma+\alpha\gamma_s\right\vert\alpha:S^1\to\mathbb{R}\right\},
\end{gather*}
and we can def\/ine a~presymplectic form $\omega_0$ on $\mathcal{M}$ by
\begin{gather*}
\omega_0(X,Y)=\int_{S^1}\det
\begin{pmatrix}
X
\\
Y
\end{pmatrix}
ds,
\qquad
 X,Y\in T_{\gamma}\mathcal{M}.
\end{gather*}
When~$X$ and~$Y$ are given by
\begin{gather}\label{Eq:2.5}
X=-\frac{1}{2}\alpha_s\gamma+\alpha\gamma_s,
\qquad
Y=-\frac{1}{2}\beta_s\gamma+\beta\gamma_s,
\qquad
\alpha,\beta: \ S^1\to\mathbb{R},
\end{gather}
a~direct calculation shows that{\samepage
\begin{gather*}
\omega_0(X,Y)=\int_{S^1}\alpha\beta_s ds,
\end{gather*}
from which we see that the kernel of $\omega_0$ at $\gamma$ is $\mathbb{R}\cdot\gamma_s$.}

It is known that the higher KdV equation~\eqref{Eq:2.4} as well as~\eqref{Eq:2.2} has an inf\/inite series of conserved
quantities $\{H_m\}_{m\in\mathbb{N}}$ given in the form of
\begin{gather*}
H_m=\int_{S^1}h_m(\kappa,\kappa_s,\kappa_{ss},\dots)ds,
\end{gather*}
where $h_m$ is a~polynomial in $\kappa$ and its derivatives up to order~$m$, for example,
\begin{gather*}
h_1=\kappa,
\qquad
h_2=\frac{1}{2}\kappa^2,
\qquad
h_3=\frac{1}{2}\kappa^3-\frac{1}{4}\kappa_s^2
\end{gather*}
(see \cite{KMGZ,L,MGK,N}).
Moreover, by using the conserved quantity,~$n$th KdV equation~\eqref{Eq:2.4} can be expressed as
\begin{gather}\label{Eq:KdVUsingVD}
\kappa_t=D_s\frac{\delta H_{n+2}}{\delta\kappa},
\end{gather}
where $\delta H_{n+2}/\delta\kappa$ is the \textit{variational derivative} of $H_{n+2}$:
\begin{gather*}
\frac{\delta H_{n+2}}{\delta\kappa}=\frac{\partial h_{n+2}}{\partial\kappa} -D_s\frac{\partial
h_{n+2}}{\partial\kappa_s}+D_s^2\frac{\partial h_{n+2}}{\partial \kappa_{ss}}-\cdots.
\end{gather*}
The expression~\eqref{Eq:KdVUsingVD} played an important role in computation in~\cite{FK3}, where we studied the higher
KdV f\/lows on the space of closed equicentroaf\/f\/ine curves as Hamiltonian systems; using the above presymplectic structure
$\omega_0$, we gave the Hamiltonian f\/lows associated with the higher KdV equations.
The paper~\cite{FK3} deals also with the geometric Miura transformation as is mentioned in Section~\ref{section5} below.

For each $n\in\mathbb{N}$, we def\/ine a~vector f\/ield $X_n$ on $\mathcal{M}$ by
\begin{gather*}
(X_n)_{\gamma}=-\frac{1}{2}\big(\Omega^{n-1}\kappa_s\big)\gamma+\big(D_s^{-1}\Omega^{n-1} \kappa_s\big)\gamma_s,
\qquad
\gamma\in\mathcal{M}.
\end{gather*}
Regarding $\{H_m\}_{m\in\mathbb{N}}$ as functions on $\mathcal{M}$ by substituting the equicentroaf\/f\/ine curvature of
$\gamma$ for $\kappa$, we have the following proposition, which is essentially due to Pinkall~\cite{P} in the case $n =
1$.
\begin{Proposition}[\cite{FK3}]
For each $n\in\mathbb{N}$, $X_n$ is a~Hamiltonian vector field for $H_n$ with respect to~$\omega_0$, i.e.,
$dH_n=\omega_0(X_n,\,\cdot\,)$ holds.
Hence $H_n$ is a~Hamiltonian function for the~$n$th KdV flow $\gamma_t=X_n$.
\end{Proposition}

Now, we def\/ine another form $\omega_1$ on $\mathcal{M}$ by
\begin{gather*}
\omega_1(X, Y)=\int_{S^1}\det
\begin{pmatrix}
X
\\
(D_s^2 + \kappa)Y
\end{pmatrix}
ds,
\qquad
X,Y\in T_{\gamma}\mathcal{M},
\end{gather*}
which is represented as
\begin{gather}\label{Eq:2.7}
\omega_1(X,Y) = \int_{S^1}\alpha\Omega\beta_s ds
\end{gather}
for $X$, $Y$ given by~\eqref{Eq:2.5}.
The following shows that $\omega_0$ and $\omega_1$ with $\{H_m\}_{m\in N}$ def\/ine a~bi-Hamiltonian structure on
$\mathcal{M}$ (cf.~\cite{M,O}).

\begin{Th}\label{Th:KdV1}
The form $\omega_1$ is a~presymplectic form on $\mathcal{M}$.
For each $n\in\mathbb{N}$, $X_n$ is a~Hamiltonian vector field for $H_{n+1}$ with respect to $\omega_1$.
\end{Th}
\begin{proof}
For two functions~$F$ and~$G$ on $\mathcal{M}$ of the form
\begin{gather}\label{Eq:FandG}
F=\int_{S^1}f(\kappa,\kappa_s,\kappa_{ss},\dots)ds,
\qquad
G=\int_{S^1}g(\kappa,\kappa_s,\kappa_{ss},\dots)ds,
\end{gather}
we set
\begin{gather*}
\{F,G\}_1=\int_{S^1}\frac{\delta F}{\delta\kappa}\Omega D_s\frac{\delta G}{\delta\kappa}ds.
\end{gather*}
Then from~\cite{KR,M}, we see that $\{\,\cdot\,,\,\cdot\,\}_1$ provides a~Poisson bracket with
\begin{gather*}
X_n = - \{H_{n + 1},\, \cdot\, \}_1.
\end{gather*}

We put $\CHVP{F} = \delta F / \delta\kappa$ and $(\HVP{F})_\gamma = - (1 / 2)(\CHVP{F})_s\gamma + \CHVP{F}\gamma_s$.
Since the dif\/ferentiation of~$F$ along a~motion $\gamma_t = X_\gamma = - (1 / 2)\alpha_s\gamma + \alpha\gamma_s$ is
given as
\begin{gather*}
XF = \frac{dF}{dt} = \int_{S^1} \frac{\delta F}{\delta\kappa}\kappa_t ds = \int_{S^1} \frac{\delta
F}{\delta\kappa}\Omega\alpha_s ds,
\end{gather*}
we have
\begin{gather*}
\omega_1(\HVP{F}, X) = \int_{S^1} \frac{\delta F}{\delta\kappa}\Omega\alpha_s ds = XF = dF(X)
\end{gather*}
and
\begin{gather*}
\omega_1(\HVP{F}, \HVP{G}) = \int_{S^1} \frac{\delta F}{\delta\kappa}\Omega D_s\frac{\delta G}{\delta\kappa} = \{F,G\}_1.
\end{gather*}
Hence $\omega_1$ is skew-symmetric and its closedness follows from the Jacobi identity for $\{\,\cdot\,,\,\cdot\,\}_1$
since for functions $F$, $G$ and $H = \int_{S^1} h(\kappa, \kappa_s, \kappa_{ss}, \dots)ds$ on $\mathcal{M}$ we have
\begin{gather*}
d\omega\big(\HVP{F}, \HVP{G}, \HVP{H}\big)=2\big(\{\{F,G\}_1,H\}_1+\{\{G,H\}_1,F\}_1+\{\{H,F\}_1,G\}_1\big)=0.
\end{gather*}
Moreover, since
\begin{gather*}
\HVP{F}G = \int_{S^1} \frac{\delta G}{\delta\kappa}\Omega D_s\frac{\delta F}{\delta\kappa} ds = \{G, F\}_1 = - \{F,
G\}_1,
\end{gather*}
we obtain $X_n = \HVP{H_{n + 1}}$ and hence
\begin{gather*}
\omega_1(X_n,\, \cdot\,) = \omega_1\big(\HVP{H_{n + 1}},\, \cdot\,\big) = dH_{n + 1}.
\end{gather*}
Therefore $X_n$ is a~Hamiltonian vector f\/ield for $H_{n+1}$ with respect to $\omega_1$.
\end{proof}

The special linear group of degree two $\mathrm{SL}(2; \mathbb{R})$ acts on $\mathcal{M}$ as $\mathcal{M} \ni \gamma
\mapsto A\gamma \in \mathcal{M}$ ($A \in \mathrm{SL}(2; \mathbb{R})$).
Two elements of $\mathcal{M}$ belong to the same orbit if and only if their equicentroaf\/f\/ine curvatures coincide.
Hence $\omega_1$ is invariant under the action of $\mathrm{SL}(2; \mathbb{R})$.
Moreover, the kernel of $\omega_1$ at $\gamma$ is the tangent space of the orbit $\mathrm{SL}(2;
\mathbb{R})\cdot\gamma$; indeed for a~one-parameter family $\gamma(\cdot, t) \in \mathcal{M}$, it follows
from~\eqref{Eq:2.2} and~\eqref{Eq:2.7} that the tangent vector~\eqref{Eq:2.1} belongs to the kernel of $\omega_1$ if and
only if $\kappa_t=0$, that is, $\kappa$ is independent of~$t$ and hence $\gamma(\cdot, t)$ is contained in an
$\mathrm{SL}(2; \mathbb{R})$-orbit.
As a~consequence, $\omega_1$ def\/ines a~symplectic form on the quotient space $\mathcal{M}/\mathrm{SL}(2; \mathbb{R})$.

We consider another action on $\mathcal{M}$ given by
\begin{gather}\label{Eq:S1-action}
\mathcal{M}\ni\gamma\mapsto\gamma(\,\cdot+\sigma)\in\mathcal{M},
\qquad
\sigma \in S^1.
\end{gather}
It is obvious that this $S^1$-action is presymplectic, that is, it leaves $\omega_1$ invariant.
Moreover, the action is Hamiltonian as we see in the proof of the following theorem.
\begin{Th}
The moment map $\mu_1$ for the $S^1$-action~\eqref{Eq:S1-action} with respect to $\omega_1$ is given by
\begin{gather}\label{Eq:MomentMapForKdV1}
\mu_1(\gamma)\left(\frac{\partial}{\partial\sigma}\right)= H_1(\gamma),
\qquad
 \gamma\in\mathcal{M} .
\end{gather}
\end{Th}

\begin{proof}
The fundamental vector f\/ield $\underline{A}$ on $\mathcal{M}$ corresponding to $\partial/\partial\sigma \in
\Lie(S^1)$ is given by $\underline{A}_\gamma=\gamma_s$ ($\gamma\in\mathcal{M}$).
For any tangent vector $\gamma_t = - (1 / 2)\alpha_s\gamma + \alpha\gamma_s$, we have
\begin{gather*}
\omega_1(\underline{A},\gamma_t) = \omega_1(\gamma_s,\gamma_t) = \int_{S^1}\Omega\alpha_sds = \int_{S^1} \kappa_t ds =
\frac{d}{dt}H_1(\gamma) = dH_1(\gamma_t),
\end{gather*}
which implies~\eqref{Eq:MomentMapForKdV1} by the def\/inition of the moment map.
\end{proof}

\begin{Rem}
Let $\Phi_n^\tau$ be the f\/low generated by $X_n$, that is, $\Phi_n^\cdot$ is a~one-parameter transformation group of
$\mathcal{M}$ such that
\begin{gather*}
\left.\frac{\partial}{\partial\tau}\right\vert_{\tau = 0}\Phi_n^\tau(\gamma) =(X_n)_{\gamma},
\qquad
\gamma\in\mathcal{M}.
\end{gather*}
As an $\mathbb{R}$-action on $\mathcal{M}$, $\Phi_n^\cdot$ is Hamiltonian with respect to $\omega_0$ and the
corresponding moment map is given by~$H_n$.
\end{Rem}

\section{Multi-Hamiltonian structures on the level sets\\ of Hamiltonians}

For a~given sequence of real numbers $C = \{c_k\}_{k \in \mathbb{N}}$, we def\/ine subsets $\mathcal{M}(C_m)$
$(m=1,2,\dots)$ of~$\mathcal{M}$ by
\begin{gather*}
\mathcal{M}(C_m) = H_1^{-1}(c_1) \cap \dots \cap H_m^{-1}(c_m).
\end{gather*}
In the following, we assume that each $\mathcal{M}(C_m)$ is not an empty set.

\begin{Lemma}
For functions $\alpha$, $\beta$ on $S^1$, if $D_s^{- 1}\Omega D_s\alpha$ is determined as a~function on $S^1$, then we have
\begin{gather}\label{Eq:3.1}
\int_{S^1} \big(D_s^{- 1}\Omega D_s\alpha\big)\cdot\beta_s ds = \int_{S^1} \alpha\Omega\beta_s ds.
\end{gather}
\end{Lemma}
\begin{proof}
Noting $\Omega D_s = (1 / 2)D_s^3 + \kappa D_s + D_s\kappa$, we can easily verify~\eqref{Eq:3.1} by integration by parts.
\end{proof}

\begin{Proposition}\label{Prop:3a}
For $\gamma \in \mathcal{M}(C_m)$ and $X = - (1 / 2)\alpha_s\gamma + \alpha\gamma_s \in T_\gamma\mathcal{M}(C_m)$,
$\Omega\alpha_s, \Omega^2\alpha_s, \dots, \linebreak[1] \Omega^{m + 1}\alpha_s$ are defined as functions on $S^1$
and $\int_{S^1} \Omega^k\alpha_s ds = 0$ for any $k = 1, 2, \dots, m$.
\end{Proposition}

\begin{proof}
We shall show the proposition by induction on~$m$.
In the case $m = 1$, $\Omega\alpha_s = (1 / 2)\alpha_{sss} + 2\kappa\alpha_s + \kappa_s\alpha$ is a~function on $S^1$
and we have
\begin{gather*}
\int_{S^1} \Omega\alpha_s ds = \int_{S^1} \kappa\alpha_s ds = \omega_0(X_1, X) = dH_1(X),
\end{gather*}
which vanishes since $X \in T_\gamma\mathcal{M}(C_1) = \mathop{\mathrm{Ker}}\nolimits (dH_1)_\gamma$; moreover, this
implies that $D_s^{- 1}\Omega\alpha_s$,
and consequently $\Omega^2\alpha_s = ((1 / 2)D_s^2 + 2\kappa+\kappa_sD_s^{- 1})\Omega\alpha_s$ are def\/ined on $S^1$.

We assume that the proposition holds for $m = l$ for some $l \geq 1$.
Then, for $X \in T_\gamma\mathcal{M}(C_{l + 1}) = T_\gamma\mathcal{M}(C_l) \cap \mathop{\mathrm{Ker}}\nolimits(dH_{l +
1})_\gamma$, by using~\eqref{Eq:3.1} we get
\begin{gather*}
\begin{split}
& \int_{S^1}\Omega^{l + 1}\alpha_s ds= \int_{S^1} \kappa\Omega^l\alpha_s ds
=\int_{S^1}\big(\big(D_s^{-1}\Omega D_s\big)^l\kappa\big)\cdot\alpha_s ds = \int_{S^1} \big(D_s^{-1}\Omega^l\kappa_s\big)\cdot\alpha_s ds
\\
&\hphantom{\int_{S^1}\Omega^{l + 1}\alpha_s ds}{}
= \omega_0(X_{l+1}, X) = dH_{l+1}(X) = 0,
\end{split}
\end{gather*}
which implies that $\Omega^{l + 2}\alpha_s$ is determined as a~function on~$S^1$ in the same way as in the case $m = 1$.
\end{proof}

From Proposition~\ref{Prop:3a}, we can def\/ine a~tensor f\/ield $\omega_{m + 1}$ of type $(0, 2)$ on $\mathcal{M}(C_m)$ by
\begin{gather*}
\omega_{m + 1}(X, Y) = \int_{S^1} \alpha\Omega^{m + 1}\beta_s ds,
\end{gather*}
which is shown to be skew-symmetric by using~\eqref{Eq:3.1}.
Furthermore, in a~similar way to the proof of Theorem~\ref{Th:KdV1}, we see that $\omega_{m + 1}$ is a~presymplectic form and $X_n$
is a~Hamiltonian vector f\/ield for the Hamiltonian function $H_{n + m + 1}$ with respect to $\omega_{m + 1}$; indeed, for
functions $F$, $G$ given by~\eqref{Eq:FandG} and for an integer~$k$, putting
\begin{gather*}
\{F,G\}_k = \int_{S^1} \frac{\delta F}{\delta\kappa}\Omega^k D_s\frac{\delta G}{\delta\kappa} ds,
\end{gather*}
we have a~family of Poisson brackets $\{\,\cdot\,,\,\cdot\,\}_k$ with
\begin{gather*}
X_n = - \{H_{n + 2 - k},\, \cdot\, \}_k.
\end{gather*}
Setting $\CHVP{F} = D_s^{- 1}\Omega^{- m}D_s(\delta F / \delta\kappa)$ and $(\HVP{F})_\gamma = - (1 /
2)(\CHVP{F})_s\gamma + \CHVP{F}\gamma_s$, we have
\begin{gather*}
\omega_{m + 1}\big(\HVP{F}, X\big) = dF(X)
\qquad \text{and}\qquad
\omega_{m + 1}\big(\HVP{F}, \HVP{G}\big) = \{F, G\}_{1 - m},
\end{gather*}
which implies that $\omega_{m + 1}$ is presymplectic.
Moreover, since
\begin{gather*}
\HVP{F}G = - \{F, G\}_{1 - m}
\end{gather*}
holds, we have $X_n = \HVP{H_{n + m + 1}}$ and
\begin{gather*}
\omega_{m + 1}(X_n,\, \cdot\,) = \omega_{m + 1}\big(\HVP{H_{n + m + 1}},\, \cdot\,\big) = dH_{n + m + 1}.
\end{gather*}
Hence $X_n$ is a~Hamiltonian vector f\/ield of $H_{n + m + 1}$ with respect to $\omega_{m + 1}$.

Besides $\omega_{m + 1}$, we have $m + 1$ more presymplectic forms on $\mathcal{M}(C_m)$ by restricting $\omega_0$,
$\omega_1$ on~$\mathcal{M}$ and $\omega_{k + 1}$'s on $\mathcal{M}(C_k)$'s for $k = 1, 2, \dots, m - 1$ to
$\mathcal{M}(C_m)$; we denote them by the same symbols.
By the discussion so far, we obtain the following theorem.

\begin{Th}\label{Th:MHSForKdV}
On $\mathcal{M}(C_m)$, for each $n \in \mathbb{N}$ and $k = 0, 1, \dots, m + 1$, $X_n$ is a~Hamiltonian vector field for
$H_{n+k}$ with respect to $\omega_k$, that is, the set $\bigl(\{H_n\}_{n \in \mathbb{N}}, \{\omega_k\}_{k = 0}^{m +
1}\bigr)$ is a~multi-Hamiltonian system on $\mathcal{M}(C_m)$ describing the higher KdV flows.
\end{Th}

As on $\mathcal{M}$, we have the following theorem for a~Hamiltonian $S^1$-action on $\mathcal{M}(C_m)$:
\begin{gather*}
\mathcal{M}(C_m)\ni\gamma\mapsto\gamma(\,\cdot+\sigma)\in\mathcal{M}(C_m),
\qquad
\sigma \in S^1.
\end{gather*}
\begin{Th}\label{Th:m+1-stMomentMap}
The moment map $\mu_{m+1}$ for the $S^1$-action on $\mathcal{M}(C_m)$ with respect to $\omega_{m+1}$ is given~by
\begin{gather*}
\mu_{m+1}(\gamma)\left(\frac{\partial}{\partial\sigma}\right)= H_{m+1}(\gamma),
\qquad
\gamma\in\mathcal{M}(C_m).
\end{gather*}
\end{Th}

\begin{Rem}
We can def\/ine $\omega_{m + 1}$ in a~manner similar to the def\/initions of $\omega_0$ and $\omega_1$.
We put a~map $\phi$ from $T_{\gamma}\mathcal{M}$ to the space of all vector f\/ields along $\gamma$ as
\begin{gather*}
\phi X = - \alpha_s\gamma,
\qquad
 X = - \frac{1}{2}\alpha_s\gamma + \alpha\gamma_s .
\end{gather*}
For any tangent vector~$X$ of $\mathcal{M}$, $(D_s^2 + \kappa)X$ has no $\gamma_s$-component and it belongs to the image
of $\phi$ if~$X$ is tangent to $\mathcal{M}(C_1)$.
Then for $X \in T_{\gamma}\mathcal{M}(C_1)$ we have
\begin{gather*}
\phi^{- 1}\big(D_s^2 + \kappa\big)X = - \frac{1}{2}(\Omega\alpha_s)\gamma + \big(D_s^{- 1}\Omega\alpha_s\big)\gamma_s
\end{gather*}
and
\begin{gather*}
\big(D_s^2 + \kappa\big)\phi^{- 1}\big(D_s^2 + \kappa\big)X = - \big(\Omega^2\alpha_s\big)\gamma.
\end{gather*}
Hence
\begin{gather*}
\int_{S^1} \det
\begin{pmatrix}
X
\\
(D_s^2 + \kappa)\phi^{- 1}(D_s^2 + \kappa)Y
\end{pmatrix}
ds = \omega_2(X, Y)
\end{gather*}
holds.
More generally, $[\phi^{- 1}(D_s^2 + \kappa)]^{m}X$ can be def\/ined for any tangent vector~$X$ of
$\mathcal{M}(C_m)$ and we obtain
\begin{gather*}
\int_{S^1} \det
\begin{pmatrix}
X
\\
(D_s^2 + \kappa)\bigl[\phi^{- 1}(D_s^2 + \kappa)\bigr]^{m}Y
\end{pmatrix}
ds = \omega_{m + 1}(X, Y)
\end{gather*}
on $\mathcal{M}(C_m)$.
We note that this formula is valid in the case $\omega_1\ (m = 0)$ and even in the case $\omega_0$ $(m=-1)$ since
\begin{gather*}
\int_{S^1} \det
\begin{pmatrix}
X
\\
\phi Y
\end{pmatrix}
ds = \omega_0(X, Y).
\end{gather*}
\end{Rem}

\section{A bi-Hamiltonian structure on the space of closed curves\\ in the Euclidean plane}\label{Section4}

We denote by $\mathbb{E}^2$ the Euclidean plane equipped with the standard inner product $\langle\, \cdot\,,\, \cdot\,\rangle$,
and we set the space $\hat{\mathcal{M}}$ of closed curves in the Euclidean plane $\mathbb{E}^2$ by
\begin{gather*}
\hat{\mathcal{M}} = \big\{\hat\gamma: S^1 \to \mathbb{E}^2 \,\big|\,\langle\hat\gamma_s(s), \hat\gamma_s(s)\rangle
\equiv 1 \big\}.
\end{gather*}
For $\hat\gamma \in \hat{\mathcal{M}}$, the curvature $\hat\kappa$ is def\/ined by $T_s = \hat\kappa N$, where $T =
\hat\gamma_s$ is the velocity vector f\/ield and~$N$ is the left-oriented unit normal vector f\/ield along $\hat\gamma$.

Let $\hat\gamma(\, \cdot\,,t)\in\hat{\mathcal{M}}$ be a~one-parameter family of closed curves in~$\mathbb{E}^2$.
Then $\hat\gamma_t$ is represented as
\begin{gather*}
\hat\gamma_t =\lambda T+\mu N,
\qquad
 \lambda,\mu: \ S^1\to\mathbb{R},
\quad
\lambda_s=\hat\kappa\mu,
\end{gather*}
and the curvature $\hat\kappa$ evolves as
\begin{gather*}
\hat\kappa_t=\mu_{ss}+\hat\kappa\lambda_s+\hat\kappa_s\lambda =\hat{\Omega}(2\mu),
\end{gather*}
where
\begin{gather*}
\hat{\Omega}=\frac{1}{2}\big(D_s^2+\hat\kappa^2+\hat\kappa_s D_s^{- 1}\hat\kappa\big)
\end{gather*}
is the \textit{recursion operator of the mKdV equation}:
\begin{gather*}
\hat\kappa_t =\hat{\Omega}\hat\kappa_s=\frac{1}{2}\hat\kappa_{sss}+\frac{3}{4} {\hat\kappa}^2\hat\kappa_s.
\end{gather*}
Hence when we choose $\mu=(1 / 2)\hat{\Omega}^{n - 1}\hat\kappa_s$, we have the~$n$th mKdV equation for $\hat\kappa$:
\begin{gather}\label{Eq:5.2}
\hat\kappa_t=\hat{\Omega}^n\hat\kappa_s.
\end{gather}

The tangent space of $\hat{\mathcal{M}}$ at $\hat\gamma\in\hat{\mathcal{M}}$ is described as
\begin{gather*}
T_{\hat\gamma}\hat{\mathcal{M}}=\{\lambda T+\mu N\, \vert\, \lambda,\mu:S^1\to\mathbb{R},\; \lambda_s= \hat\kappa\mu\},
\end{gather*}
and we can def\/ine a~presymplectic form $\hat\omega_0$ on $\hat{\mathcal{M}}$ by
\begin{gather*}
\hat\omega_0(X,Y)=\int_{S^1}\langle D_sX,Y\rangle ds,
\qquad
X,Y\in T_{\hat\gamma}\hat{\mathcal{M}}.
\end{gather*}
When~$X$ and~$Y$ are given by
\begin{gather}\label{Eq:5.3}
X=\lambda T+\mu N,
\qquad
Y=\tilde\lambda T+\tilde\mu N,
\qquad
\lambda,\mu,\tilde\lambda,\tilde\mu: \ S^1\to\mathbb{R},
\end{gather}
we have
\begin{gather*}
\hat\omega_0(X,Y)=\int_{S^1}(\hat\kappa\lambda+\mu_s)\tilde{\mu}ds,
\end{gather*}
and we see that the kernel of $\hat\omega_0$ at $\hat\gamma$ is $\mathbb{R}\cdot\hat\gamma_s$.

As in the case of the higher KdV equation~\eqref{Eq:2.4}, the~$n$th mKdV equation~\eqref{Eq:5.2} can be written as
\begin{gather*}
\hat\kappa_t=D_s\frac{\delta\hat H_{n+2}}{\delta\hat\kappa}
\end{gather*}
for an inf\/inite series of conserved quantities $\{\hat H_m\}_{m\in\mathbb{N}}$ expressed in the form of
\begin{gather*}
\hat H_m=\int_{S^1}\hat h_m(\hat\kappa,\hat\kappa_s,\hat\kappa_{ss},\dots)ds,
\end{gather*}
where $\hat h_m$ is a~polynomial in $\hat\kappa$ and its derivatives up to order~$m$, for example,
\begin{gather*}
\hat h_1=\frac{1}{4}\hat\kappa^2,
\qquad
\hat h_2=\frac{1}{32}\hat\kappa^4-\frac{1}{8}\hat\kappa_s^2,
\qquad
\hat h_3=\frac{1}{128}\hat\kappa^6-\frac{5}{32}\hat\kappa^2\hat\kappa_s^2 +\frac{1}{16}\hat\kappa_{ss}^2.
\end{gather*}

For each $n\in\mathbb{N}$, we def\/ine a~vector f\/ield $\hat X_n$ on $\hat{\mathcal{M}}$ by
\begin{gather*}
\big(\hat X_n\big)_{\hat\gamma}
=\frac{1}{2}\big(D_s^{-1}\big(\hat\kappa{\hat{\Omega}}^{n-1}\hat\kappa_s\big)\big)T+\frac{1}{2}\big({\hat{\Omega}}^{n-1}\hat\kappa_s\big)N,
\qquad
 \hat\gamma \in \hat{\mathcal{M}} ,
\end{gather*}
then we have the following.
\begin{Proposition}[\cite{FK3}]
For each $n\in\mathbb{N}$, $\hat X_n$ is a~Hamiltonian vector field for $\hat H_n$ with respect to~$\hat\omega_0$.
Hence $\hat H_n$ is a~Hamiltonian function for the~$n$th mKdV flow $\hat\gamma_t= \hat X_n$.
\end{Proposition}
In addition, we def\/ine another form $\hat\omega_1$ on $\hat{\mathcal{M}}$ by
\begin{gather*}
\hat\omega_1(X, Y)=\int_{S^1}\big\langle D_sX,D_s^2Y\big\rangle ds,
\qquad
X,Y\in T_{\hat\gamma}\hat{\mathcal{M}},
\end{gather*}
which is represented as
\begin{gather*}
\hat\omega_1(X,Y)=\int_{S^1}(\hat\kappa\lambda+\mu_s)\hat{\Omega}\tilde{\mu}ds
\end{gather*}
for $X$, $Y$ given by~\eqref{Eq:5.3}.
The following theorem is proved in a~similar way to the proof of Theorem~\ref{Th:KdV1}.
\begin{Th}
The form $\hat\omega_1$ is a~presymplectic form on $\hat{\mathcal{M}}$.
For each $n\in\mathbb{N}$, $\hat X_n$ is a~Hamiltonian vector field for $\hat H_{n+1}$ with respect to $\hat\omega_1$.
\end{Th}

Note that the Euclidean motion group $E(2) = O(2) \ltimes \mathbb{R}^2$ of $\mathbb{E}^2$ acts on $\hat{\mathcal{M}}$.
It is easily verif\/ied that $\hat{\omega}_1$ is invariant under the $E(2)$-action and the kernel of $\hat{\omega}_1$ at
$T_{\hat\gamma}\hat{\mathcal{M}}$ contains the tangent space of the orbit.
Hence $\omega_1$ determines a~presymplectic form on $\hat{\mathcal{M}}/E(2)$.

As well as on $(\mathcal{M}, \omega_1)$, $S^1$ acts on $\hat{\mathcal{M}}$ leaving $\hat{\omega}_1$ invariant and the
following theorem holds.
\begin{Th}
The moment map $\hat\mu_1$ for the $S^1$-action on $\hat{\mathcal{M}}$ with respect to $\hat{\omega}_1$ is given by
\begin{gather*}
\hat\mu_1(\hat\gamma)\left(\frac{\partial}{\partial\sigma}\right) =\hat{H}_1(\hat\gamma),
\qquad
 \hat\gamma\in\hat{\mathcal{M}} .
\end{gather*}
\end{Th}

\section{The geometric Miura transformation and multi-Hamiltonian structures on spaces of closed curves in the Euclidean plane}\label{section5}

First, we brief\/ly review the geometric Miura transformation which relates the Hamiltonian structures on $\mathcal{M}$
and on $\hat{\mathcal{M}}$ (see~\cite{FK3} for more details).
We consider the complexif\/ication of $\mathcal{M}$:
\begin{gather*}
\mathcal{M}^\mathbb{C} = \left\{\gamma: S^1 \to \mathbb{C}^2 \setminus \{0\}\left\vert \det
\begin{pmatrix}
\gamma
\\
\gamma_s
\end{pmatrix}
= 1\right.\right\}.
\end{gather*}
We determine the curvature of $\gamma \in \mathcal{M}^\mathbb{C}$, (complex) presymplectic forms on
$\mathcal{M}^\mathbb{C}$, etc.\ by the same formulas as in the case of $\mathcal{M}$, hence we use the same symbols
$\kappa, \omega_0, \omega_1, \dots$ to denote them.

By identifying the range $\mathbb{E}^2$ of $\hat{\gamma} \in \hat{\mathcal{M}}$ with a~complex plane $\mathbb{C}$, we
def\/ine the {\it geometric Miura transformation} $\Phi: \hat{\mathcal{M}} \to \mathcal{M}^\mathbb{C}$ by
\begin{gather*}
\Phi(\hat{\gamma}) = (- \hat{\gamma}_s)^{-\frac{1}{2}}\left(\hat{\gamma}, 1\right),
\qquad
 \hat{\gamma} \in \hat{\mathcal{M}}.
\end{gather*}
The curvature $\kappa$ of $\Phi(\hat{\gamma})$ is related with the curvature $\hat{\kappa}$ of $\hat{\gamma}$ by the
Miura transformation:
\begin{gather}\label{Eq:6:1}
\kappa=\frac{\sqrt{- 1}}{2}\hat{\kappa}_s+\frac{1}{4}\hat{\kappa}^2.
\end{gather}
Moreover, we have the following.
\begin{Proposition}[\cite{FK3}]
For each $n\in\mathbb{N}$, $\Phi_*\hat{X}_n = X_n$ holds and the Hamiltonian system $(\hat{\omega}_0,
\hat{H}_n)$ on $\hat{\mathcal{M}}$ coincides with the pullback of $(\omega_0, H_n)$ on $\mathcal{M}^\mathbb{C}$ by~$\Phi$:
\begin{gather}\label{Eq:6:2}
\hat{\omega}_0 = \Phi^*\omega_0,
\qquad
\hat{H}_n = \Phi^*H_n.
\end{gather}
\end{Proposition}

For a~sequence of real numbers $C = \{c_k\}_{k \in \mathbb{N}}$, the second equation of~\eqref{Eq:6:2} implies that
\begin{gather*}
\hat{\mathcal{M}}(C_m) = \hat{H}_1^{-1}(c_1) \cap \dots \cap \hat{H}_m^{-1}(c_m) = \Phi^{-
1}\bigl(\mathcal{M}^\mathbb{C}(C_m)\bigr).
\end{gather*}
Therefore, $\Phi$ gives a~map from $\hat{\mathcal{M}}(C_m)$ to $\mathcal{M}^\mathbb{C}(C_m)$ and we have a~presymplectic
form $\hat{\omega}_{m + 1} = \Phi^*\omega_{m + 1}$ on $\hat{\mathcal{M}}(C_m)$.
Under these settings the following theorems are directly deduced from Theorems~\ref{Th:MHSForKdV} and~\ref{Th:m+1-stMomentMap}.

\begin{Th}
On $\hat{\mathcal{M}}(C_m)$, for each $n \in \mathbb{N}$ and $k = 0, 1, \dots, m + 1$, $\hat{X}_n$ is a~Hamiltonian
vector field for $\hat{H}_{n + k}$ with respect to $\hat{\omega}_k$, that is, the set $\bigl(\{\hat{H}_n\}_{n \in
\mathbb{N}}, \{\hat{\omega}_k\}_{k = 0}^{m + 1}\bigr)$ is a~multi-Hamiltonian system on $\hat{\mathcal{M}}(C_m)$
describing the higher modified KdV flows.
\end{Th}

\begin{Th}
The moment map $\hat\mu_{m+1}$ for the $S^1$-action on $\hat{\mathcal{M}}(C_m)$ with respect to $\hat{\omega}_{m+1}$ is
given by
\begin{gather*}
\hat\mu_{m+1}(\hat\gamma)\left(\frac{\partial}{\partial\sigma}\right)=\hat H_{m+1} (\hat\gamma),
\qquad
 \hat\gamma\in\hat{\mathcal{M}}(C_m) .
\end{gather*}
\end{Th}

\begin{Rem}
The symplectic form $\omega_{m + 1}$ can be represented as
\begin{gather}\label{Eq:Expliciteomegahat}
\hat{\omega}_{m + 1}(X, Y) = \int_{S^1} (\hat{\kappa}\lambda + \mu_s)\hat{\Omega}^{m + 1}\tilde{\mu} ds,
\end{gather}
where~$X$ and~$Y$ are tangent vectors on $\hat{\mathcal{M}}(C_m)$ given by~\eqref{Eq:5.3}.
In fact, when $\kappa$ and $\hat{\kappa}$ are related by~\eqref{Eq:6:1}, a~direct calculation shows an identity
\begin{gather*}
\left(\sqrt{- 1} D_s + \hat{\kappa}\right)\hat{\Omega} = \Omega\left(\sqrt{- 1} D_s + \hat{\kappa}\right);
\end{gather*}
thus we have
\begin{gather*}
\hat{\omega}_{m + 1}(X, Y)= \omega_{m+1}\bigl(\Phi_*X,\Phi_*Y\bigr)
= \int_{S^1}\big(\lambda+\sqrt{- 1}\mu\big)\Omega^{m+1}\big(\tilde{\lambda}+\sqrt{- 1}\tilde{\mu}\big)_s ds
\\
\phantom{\hat{\omega}_{m + 1}(X, Y)}{}
=\int_{S^1} \big(\lambda + \sqrt{- 1}\mu\big)\Omega^{m + 1}\big(\sqrt{- 1} D_s + \hat{\kappa}\big)\tilde{\mu} ds
\\
\phantom{\hat{\omega}_{m + 1}(X, Y)}{}
= \int_{S^1} \big(\lambda + \sqrt{- 1}\mu\big)\big(\sqrt{- 1} D_s + \hat{\kappa}\big)\hat{\Omega}^{m + 1}\tilde{\mu} ds
\\
\phantom{\hat{\omega}_{m + 1}(X, Y)}{}
= \int_{S^1} \bigl[\big({-} \sqrt{- 1} D_s + \hat{\kappa}\big)\big(\lambda + \sqrt{- 1}\mu\big)\bigr]\cdot\hat{\Omega}^{m +
1}\tilde{\mu} ds
\\
\phantom{\hat{\omega}_{m + 1}(X, Y)}
=\int_{S^1}(\hat\kappa\lambda+\mu_s)\hat{\Omega}^{m + 1}\tilde{\mu}ds.
\end{gather*}
We note that~\eqref{Eq:Expliciteomegahat} implies $\hat{\omega}_{m + 1}$ is a~real form, though $\omega_{m + 1}$ on
$\mathcal{M}^\mathbb{C}(C_m)$ is complex.
\end{Rem}

\subsection*{Acknowledgements}

The authors would like to thank the referees' kind and important comments and advice.
The f\/irst named author is partly supported by JSPS KAKENHI Grant Number 22540070 and by the Kansai University
Grant-in-Aid for progress of research in graduate course, 2013.
The second named author is partly supported by JSPS KAKENHI Grant Number 22540107.

\pdfbookmark[1]{References}{ref}
\LastPageEnding


\begin{thebibliography}{99}
\footnotesize \itemsep=0pt

\bibitem{A1}
Anco S.C., Bi-{H}amiltonian operators, integrable f\/lows of curves using moving
  frames and geometric map equations, \href{http://dx.doi.org/10.1088/0305-4470/39/9/005}{\textit{J.~Phys.~A: Math. Gen.}}
  \textbf{39} (2006), 2043--2072, \href{http://arxiv.org/abs/nlin.SI/0512051}{nlin.SI/0512051}.

\bibitem{A2}
Anco S.C., Hamiltonian f\/lows of curves in {$G/{\rm SO}(N)$} and vector soliton
  equations of m{K}d{V} and sine-{G}ordon type, \href{http://dx.doi.org/10.3842/SIGMA.2006.044}{\textit{SIGMA}} \textbf{2}
  (2006), 044, 18~pages, \href{http://arxiv.org/abs/nlin.SI/0512046}{nlin.SI/0512046}.

\bibitem{A3}
Anco S.C., Group-invariant soliton equations and bi-{H}amiltonian geometric
  curve f\/lows in {R}iemannian symmetric spaces, \href{http://dx.doi.org/10.1016/j.geomphys.2007.09.005}{\textit{J.~Geom. Phys.}}
  \textbf{58} (2008), 1--37, \href{http://arxiv.org/abs/nlin.SI/0703041}{nlin.SI/0703041}.

\bibitem{A4}
Anco S.C., Hamiltonian curve f\/lows in {L}ie groups {$G\subset {\rm U}(N)$} and
  vector {NLS}, m{K}d{V}, sine-{G}ordon soliton equations, in Symmetries and
  Overdetermined Systems of Partial Dif\/ferential Equations, \href{http://dx.doi.org/10.1007/978-0-387-73831-4_10}{\textit{IMA Vol.
  Math. Appl.}}, Vol.~144, Springer, New York, 2008, 223--250,
  \href{http://arxiv.org/abs/nlin.SI/0610075}{nlin.SI/0610075}.

\bibitem{AA1}
Anco S.C., Asadi E., Quaternionic soliton equations from {H}amiltonian curve
  f\/lows in {${\mathbb{HP}}^n$}, \href{http://dx.doi.org/10.1088/1751-8113/42/48/485201}{\textit{J.~Phys.~A: Math. Theor.}} \textbf{42}
  (2009), 485201, 25~pages, \href{http://arxiv.org/abs/0905.4215}{arXiv:0905.4215}.

\bibitem{AA2}
Anco S.C., Asadi E., Symplectically invariant soliton equations from
  non-stretching geometric curve f\/lows, \href{http://dx.doi.org/10.1088/1751-8113/45/47/475207}{\textit{J.~Phys.~A: Math. Theor.}}
  \textbf{45} (2012), 475207, 37~pages, \href{http://arxiv.org/abs/1206.4040}{arXiv:1206.4040}.

\bibitem{AM}
Anco S.C., Myrzakulov R., Integrable generalizations of {S}chr\"odinger maps
  and {H}eisenberg spin models from {H}amiltonian f\/lows of curves and surfaces,
  \href{http://dx.doi.org/10.1016/j.geomphys.2010.05.013}{\textit{J.~Geom. Phys.}} \textbf{60} (2010), 1576--1603, \href{http://arxiv.org/abs/0806.1360}{arXiv:0806.1360}.

\bibitem{AV}
Anco S.C., Vacaru S.I., Curve f\/lows in {L}agrange--{F}insler geometry,
  bi-{H}amiltonian structures and solitons, \href{http://dx.doi.org/10.1016/j.geomphys.2008.10.006}{\textit{J.~Geom. Phys.}} \textbf{59}
  (2009), 79--103, \href{http://arxiv.org/abs/math-ph/0609070}{math-ph/0609070}.

\bibitem{CIB}
Calini A., Ivey T., Mar{\'{\i}}-Bef\/fa G., Remarks on {K}d{V}-type f\/lows on
  star-shaped curves, \href{http://dx.doi.org/10.1016/j.physd.2009.01.007}{\textit{Phys.~D}} \textbf{238} (2009), 788--797,
  \href{http://arxiv.org/abs/0808.3593}{arXiv:0808.3593}.

\bibitem{CQ1}
Chou K.-S., Qu C., The {K}d{V} equation and motion of plane curves,
  \href{http://dx.doi.org/10.1143/JPSJ.70.1912}{\textit{J.~Phys. Soc. Japan}} \textbf{70} (2001), 1912--1916.

\bibitem{CQ2}
Chou K.-S., Qu C., Integrable equations arising from motions of plane curves,
  \href{http://dx.doi.org/10.1016/S0167-2789(01)00364-5}{\textit{Phys.~D}} \textbf{162} (2002), 9--33.

\bibitem{CQ3}
Chou K.-S., Qu C., Integrable motions of space curves in af\/f\/ine geometry,
  \href{http://dx.doi.org/10.1016/S0960-0779(01)00179-5}{\textit{Chaos Solitons Fractals}} \textbf{14} (2002), 29--44.

\bibitem{FK1}
Fujioka A., Kurose T., Motions of curves in the complex hyperbola and the
  {B}urgers hierarchy, \textit{Osaka~J. Math.} \textbf{45} (2008), 1057--1065.

\bibitem{FK2}
Fujioka A., Kurose T., Geometry of the space of closed curves in the complex
  hyperbola, \href{http://dx.doi.org/10.2206/kyushujm.63.161}{\textit{Kyushu~J. Math.}} \textbf{63} (2009), 161--165.

\bibitem{FK3}
Fujioka A., Kurose T., Hamiltonian formalism for the higher {K}d{V} f\/lows on
  the space of closed complex equicentroaf\/f\/ine curves, \href{http://dx.doi.org/10.1142/S0219887810003987}{\textit{Int.~J. Geom.
  Methods Mod. Phys.}} \textbf{7} (2010), 165--175.

\bibitem{H}
Hasimoto H., A soliton on a vortex f\/ilament, \href{http://dx.doi.org/10.1017/S0022112072002307}{\textit{J.~Fluid Mech.}}
  \textbf{51} (1972), 477--485.

\bibitem{KMGZ}
Kruskal M.D., Miura R.M., Gardner C.S., Zabusky N.J., Korteweg--de {V}ries
  equation and generalizations. {V}.~{U}niqueness and nonexistence of
  polynomial conservation laws, \href{http://dx.doi.org/10.1063/1.1665232}{\textit{J.~Math. Phys.}} \textbf{11} (1970),
  952--960.

\bibitem{KR}
Kulish P.P., Reiman A.G., Hierarchy of symplectic forms for the {S}chr\"odinger
  equation and for the {D}irac equation on a line, \href{http://dx.doi.org/10.1007/BF01375613}{\textit{J.~Sov. Math.}}
  \textbf{22} (1983), 1627--1637.

\bibitem{Lam}
Lamb Jr. G.L., Solitons and the motion of helical curves, \href{http://dx.doi.org/10.1103/PhysRevLett.37.235}{\textit{Phys. Rev.
  Lett.}} \textbf{37} (1976), 235--237.

\bibitem{L}
Lax P.D., Integrals of nonlinear equations of evolution and solitary waves,
  \href{http://dx.doi.org/10.1002/cpa.3160210503}{\textit{Comm. Pure Appl. Math.}} \textbf{21} (1968), 467--490.

\bibitem{LQZ}
Liu Y., Qu C., Zhang Y., Stability of periodic peakons for the modif\/ied
  {$\mu$}-{C}amassa--{H}olm equation, \href{http://dx.doi.org/10.1016/j.physd.2013.02.001}{\textit{Phys.~D}} \textbf{250} (2013),
  66--74.

\bibitem{M}
Magri F., A simple model of the integrable {H}amiltonian equation,
  \href{http://dx.doi.org/10.1063/1.523777}{\textit{J.~Math. Phys.}} \textbf{19} (1978), 1156--1162.

\bibitem{B}
Mar{\'{\i}}~Bef\/fa G., Geometric realizations of bi-{H}amiltonian completely
  integrable systems, \href{http://dx.doi.org/10.3842/SIGMA.2008.034}{\textit{SIGMA}} \textbf{4} (2008), 034, 23~pages,
  \href{http://arxiv.org/abs/0803.3866}{arXiv:0803.3866}.

\bibitem{BSW}
Mar{\'{\i}}~Bef\/fa G., Sanders J.A., Wang J.P., Integrable systems in
  three-dimensional {R}iemannian geometry, \href{http://dx.doi.org/10.1007/s00332-001-0472-y}{\textit{J.~Nonlinear Sci.}}
  \textbf{12} (2002), 143--167.

\bibitem{MGK}
Miura R.M., Gardner C.S., Kruskal M.D., Korteweg--de {V}ries equation and
  generalizations. {II}.~{E}xistence of conservation laws and constants of
  motion, \href{http://dx.doi.org/10.1063/1.1664701}{\textit{J.~Math. Phys.}} \textbf{9} (1968), 1204--1209.

\bibitem{N}
Newell A.C., Solitons in mathematics and physics, \href{http://dx.doi.org/10.1137/1.9781611970227}{\textit{CBMS-NSF Regional
  Conference Series in Applied Ma\-the\-matics}}, Vol.~48, Society for Industrial
  and Applied Mathematics (SIAM), Philadelphia, PA, 1985.

\bibitem{O}
Olver P.J., Applications of {L}ie groups to dif\/ferential equations,
  \href{http://dx.doi.org/10.1007/978-1-4684-0274-2}{\textit{Graduate Texts in Mathematics}}, Vol.~107, Springer-Verlag, New York,
  1986.

\bibitem{P}
Pinkall U., Hamiltonian f\/lows on the space of star-shaped curves,
  \href{http://dx.doi.org/10.1007/BF03322836}{\textit{Results Math.}} \textbf{27} (1995), 328--332.

\bibitem{RS}
Rogers C., Schief W.K., B\"acklund and {D}arboux transformations. Geometry and
  modern applications in soliton theory, \href{http://dx.doi.org/10.1017/CBO9780511606359}{\textit{Cambridge Texts in Applied
  Mathematics}}, Cambridge University Press, Cambridge, 2002.

\bibitem{SW}
Sanders J.A., Wang J.P., Integrable systems in {$n$}-dimensional {R}iemannian
  geometry, \textit{Mosc. Math.~J.} \textbf{3} (2003), 1369--1393,
  \href{http://arxiv.org/abs/math.AP/0301212}{math.AP/0301212}.

\bibitem{SB}
Squires S.A., Mar{\'{\i}}~Bef\/fa G., Integrable systems associated to curves in f\/lat
  {G}alilean and {M}inkowski spaces, \href{http://dx.doi.org/10.1080/00036810903397487}{\textit{Appl. Anal.}} \textbf{89} (2010),
  571--592.

\bibitem{TT}
Terng C.-L., Thorbergsson G., Completely integrable curve f\/lows on adjoint
  orbits, \href{http://dx.doi.org/10.1007/BF03322713}{\textit{Results Math.}} \textbf{40} (2001), 286--309,
  \href{http://arxiv.org/abs/math.DG/0108154}{math.DG/0108154}.

\end{thebibliography}
\end{document}